\documentclass[reqno, 12pt]{amsart}

\usepackage{amssymb}
\usepackage{epsfig}
\usepackage{epsf}
\usepackage{float}
\usepackage{a4wide}

\numberwithin{equation}{section}

\renewcommand\thefigure{\thesection.\@arabic\c@figure}
\renewcommand\thetable{\thesection.\@arabic\c@table}       

\newtheorem{thm}{Theorem}[section]
\newtheorem{lemma}[thm]{Lemma}
\newtheorem{corol}[thm]{Corollary}
\newtheorem{propos}[thm]{Proposition}

\def\reff#1{(\ref{#1})}
\begin{document}

\def\E{{\mathbb E}}
\def\P{{\mathbb P}}
\def\R{{\mathbb R}}
\def\Z{{\mathbb Z}}
\def\N{{\mathbb N}}
\def\cF{{\mathcal F}}
\def\cX{{\mathcal X}}
\def\cY{{\mathcal Y}}
\def\cN{{\mathcal N}}
\def\onu{{\overline \nu}}
\def\omu{{\overline \mu}}
\def\opi{{\overline \pi}}
\def\teta{{\widetilde\eta}}
\def\tzeta{{\widetilde\zeta}}
\def\tsigma{{\widetilde\sigma}}
\def\ttheta{{\widetilde\theta}}
\def\txi{{\widetilde\xi}}
\def\tgamma{{\widetilde\gamma}}
\def\hsigma{{\widehat\sigma}}
\def\hxi{{\widehat\xi}}
\def\hzeta{{\widehat\zeta}}
\def\hgamma{{\widehat\gamma}}
\def\G{{\mathcal T}}
\def\T{{\cal T}}
\def\I{{\cal I}}
\def\TT{\bar{{\cal T}}}
\def\II{\bar{{\cal I}}}
\def\C{{\C}}
\def\C{{\cal D}}
\def\n{{\bf n}}
\def\m{{\bf m}}
\def\b{{\bf b}}
\def\Var{{\hbox{Var}}}
\def\Cov{{\hbox{Cov}}}
\def\uR{{\underline R}}
\def\oR{{\overline R}}
\def\urho{{\underline \rho}}
\def\orho{{\overline \rho}}

\def\sqr{\vcenter{
         \hrule height.1mm
         \hbox{\vrule width.1mm height2.2mm\kern2.18mm\vrule width.1mm}
         \hrule height.1mm}}                  % This is a slimmer sqr.
\def\square{\ifmmode\sqr\else{$\sqr$}\fi}
\def\one{{\bf 1}\hskip-.5mm}
\def\limn{\lim_{N\to\infty}}
\def\given{\ \vert \ }
\def\ze{{\zeta}}
\def\be{{\beta}}
\def\la{{\lambda}}
\def\ga{{\gamma}}
\def\a{{\alpha}}
\def\th{{\theta}}
\def\proof{\noindent{\bf Proof. }}
\def\A{{\bf A}}
\def\B{{\bf B}}
\def\C{{\bf C}}
\def\D{{\bf D}}
\def\MM{{\bf m}}
\def\w{\bar{w}}
\def\lnt{{\Lambda^N}}
\def\dlnt{\delta\Lambda^N_t}
\def\lno{\Lambda^N_0}
\def\dlno{\delta\Lambda^N_0}

\title{Law of large numbers for the asymmetric simple exclusion process}
\date{}

\author{E. Andjel, P. A. Ferrari, A. Siqueira}

\maketitle

\paragraph{\bf Abstract} 
We consider simple exclusion processes on $\Z$ for which the underlying
random walk has a finite first moment and a non-zero mean and whose initial
distributions are product measures with different densities to the left and
to the right of the origin. We prove a strong law of large numbers for the
number of particles present at time $t$ in an interval growing linearly with
$t$.



\section{Introduction}
The simple exclusion process is a well studied interacting particle system.
An excellent introduction to the subject is Chapter 8 of Liggett (1985). We
study the one dimensional case with state space ${\mathcal X}:=\{0,1\}^\Z$ and
associated transition probability matrix $p(x,y)$ satisfying
\begin{enumerate}
\item $p(x,y)=p(0,y-x)$ $\forall \ x,y\in \Z,$
\item $M:=\sum_{x\in\Z} \vert x \vert p(0,x)<\infty$
%\item $\sum_{x\in\Z} x\,p(0,x)=\alpha >0.$
\end{enumerate}
The mean jump $\alpha:=\sum_{x\in\Z} x\,p(0,x)$ is called the \emph{drift}. 
 For an arbitrary initial
distribution $\mu$, we denote by $\P_{\mu}$ the probability measure on the
space of trajectories of the process associated to $\mu$.
The initial distribution will be a product measure $ \nu ^{\lambda,\rho}$ with
marginals $\nu ^{\lambda,\rho}(\{\eta : \eta(x)=1\})=\rho $ if $x> 0$ and $\nu
^{\lambda,\rho}(\{\eta : \eta(x)=1\})=\lambda$ if $x\le 0$. If $\lambda=\rho$,
this measure will also be denoted by $\nu ^{\rho}$. We recall that the
measures $\nu ^{\rho}$ ($0\leq \rho\leq 1$) are invariant for the process.
Let $\eta_t(x)$ indicate the presence/absence of a particle at site $x\in\Z$
at time $t\ge 0$. Our main result is: 

\begin{thm} 
\label{1.1}
  Let $0\leq \rho \leq \lambda \leq 1$ and let $u<v$ be real numbers. Then,
\[ \P_{\nu^{\lambda,\rho}} \Bigl( \lim_{t\to\infty} \frac {1}{t} 
\sum _{ut\le x\le vt}\eta_t (x)=   \int_u^v f(s)ds           \Bigr)\;=\;1\]
where 
\begin{equation}
  \label{x10}
  f(u)\;:=\;\left\{
\begin{array}{lllllllll}
\lambda & \hbox{if }& u<\alpha(1-2\lambda)  &  &  &  &  &  &  \\ 
\frac{1}{2}(1-\frac{u}{\alpha} ) & \hbox{if } & \alpha(1-2\lambda)\leq u\leq
\alpha (1-2\rho) & &  &  &  &  &  \\ 
\rho & \hbox{if}&  \alpha(1-2\rho)<u &  &  &  &  &  & 
\end{array}
\right. , 
\end{equation}
if $\alpha>0$ and
\begin{equation}
  \label{x11}
  f(u)\;:=\;\left\{
\begin{array}{lllllllll}
\lambda & \hbox{if }& u<\alpha(1-\lambda-\rho)  &  &  &  &  &  &  \\ 
\rho & \hbox{if}&  \alpha(1-\lambda-\rho)\leq u &  &  &  &  &  & 
\end{array}
\right. ,
\end{equation}
if $\alpha\leq 0$.
\end{thm}

Theorem \ref{1.1} means that the number of particles in the interval $[ut,vt]$
at time $t$ satisfies a strong law of large numbers as $t$ goes to infinity.
The function $f$ is the entropic solution at time 1 of the Burgers equation
associated to the exclusion process. The convergence in probability of $ \frac
{1}{t} \sum_{ut\le x\le vt}\eta_t (x)$ is proved by Rezakhanlou (1991) for a
large class of initial product measures when $p(x,y)=0$ whenever $\vert x-y
\vert$ is larger than an arbitrary constant. In the nearest neighbors totally
asymmetric case ($p(x,x+1)= 1$ and $p(x,y)=0$ otherwise) the almost sure
convergence was proven by Rost (1981) using the subadditive ergodic theorem
and by Seppalainen (1998) using a variational approach. Benassi and Fouque
(1987) stated the almost sure convergence under the nearest-neighbors
assumption $p(x,y)=0$ if $\vert x-y \vert>1$.  Unfortunately their proof
contains a mistake, as will be explained later.  Nevertheless the main ideas
of their proof, the use of a subaditive ergodic theorem and the introduction
of various classes of particles, are valuable and essential in this paper. In
those references the reader can find the relation between the exclusion
process and the Burgers equation.

\section{Graphical construction and coupling}

\paragraph{\bf Graphical construction} 
It is convenient to perform a Harris graphical construction of the process
(see Harris 1978).  Let ${\mathcal N}=((N_t({x,y}),t\ge 0)\,:\, x,y\in\Z)$ be
a family of independent Poisson processes such that the rate of the process
indexed by $(x,y)$ is $p(x,y)$. These processes are constructed on a
probability space denoted $(\Omega, {\mathcal A},\P)$. We denote $\E$ the
expectation with respect to $\P$. The rate of $\sum_{y\in\Z} N_t({x,y})$ is
$1$ for all $x$ and the rates of $\sum_{y<x} \sum_{z\geq x}N_t(y,z)$ and
$\sum_{z<x}\sum_{y\geq x}N_t(y,z)$ are bounded by $M$ for all $x$.

The process $\eta_t$ is now constructed as follows: a particle at $x\in \Z$
waits until $\sum_y N_t({x,y})$ jumps. If this jump is due to a jump of $
N_t({x,z})$ then the particle either remains at $x$ or jumps to $z$ according
to whether or not $z$ is occupied by another particle.  We will now show that
on a set of probability $1$, the trajectories of the process $\eta_t$ are well
defined for all initial configurations.  We ignore the set of probability $0$
on which two Poisson processes have a simultaneous jump.  Fix $z\in \Z$ and
$k\in \N$. The number of crossings in the time
interval $[0,k]$ of an arbitrary site $i$,  
\begin{equation}
  \label{y7}
  \sum_{x<i} \sum_{y\geq i} (N_k(x,y)+ N_k(y,x)),
\end{equation}
is a Poisson random variable with mean
\begin{equation}
  \label{y8}
  k\,\sum_{x<i} \sum_{y\geq i} (p(x,y)+p(y,x))
  \; =\; k\,\sum_x |x| p(0,x) \;=\;kM\;<\;\infty \ ,  
\end{equation}
where the inequality follows from our hypothesis on $p(x,y)$. In particular
there is a positive probability of no crossings: $\P(\sum_{x<i} \sum_{y\geq i}
(N_k(x,y)+ N_k(y,x))=0)>0$. Then, by the Ergodic Theorem, with probability $1$
there exist random integers $i$ and $j$ such that $i<z<j$ and
\[
\sum_{x<i} \sum_{y\geq i} (N_k(x,y)+ N_k(y,x))
\;=\; \sum_{x<j} \sum_{y\geq
  j} (N_k(y,z)+ N_k(z,y))\;=\;0
\]
This means that, on a set $\Omega_{z,k}$ of probability $1$, up to time $k$ no
particles have crossed (or atempted to cross) the boundaries of the interval
$[i,j-1]$. Hence the mouvement of the particles which at time $0$ were in that
interval depend up to time $k$ only on a finite number of Poisson processes
and are therefore properly determined.  Repeating the argument for each $z\in
\Z$ and each $k\in \N $ we see that on the set $\cap_{z,k}\Omega_{x,k}$,
$\eta_t(z)$ is well defined, as a function of $\mathcal N$ for all $t\geq 0$,
$z\in \Z$ and $\eta \in \mathcal X $. When the dependence on $ \mathcal N$
needs to be emphasized we will write $\eta_t(\mathcal N)$.  It is convenient
to start the process at different times, with different initial
configurations; so we also use the notation $\eta^{\eta,s}_t$ to denote the
configuration at time $t$ for a process that at time $s$ started with the
initial configuration $\eta$; when $s=0$ we just denote $\eta^\eta_t$. In
particular we have
\begin{equation}
  \label{p17}
  \eta^{\eta,s}_t({\mathcal N}) = \eta^{\eta,0}_{t-s}(\tau_s{\mathcal N}) 
\end{equation}
where $\tau_s{\mathcal N}=((N_t({x,y})-N_s({x,y}),t\ge s)\,:\, x,y\in\Z)$.
Denote 
\begin{equation}
  \label{w2}
  \P_\mu
  (\eta_t\in A) :=  \int \mu(d\eta)\P(\eta^\eta_t\in A)\,;\quad 
  \E_\mu F(\eta_t) := \int \mu(d\eta) \E(F(\eta^\eta_t))
\end{equation}
for measurable sets $A\subset \cX$ and continuous functions $F:\cX\to\R$.

\paragraph{\bf Two types of particles (coupling)} 
The graphical construction allows the simultaneous construction (with the same
Poisson processes) of various versions of the process starting with different
initial configurations. Let $\sigma$ and $\theta$ two initial configurations and
consider the \emph{coupled} process 
\begin{equation}
  \label{p38}
  (\sigma_t,\theta_t):=(\eta^\sigma_t({\mathcal
  N}),\eta^\theta_t({\mathcal N}))
\end{equation}
This is a Markov process in ${\cX}^2$.  We can continue using $\P$ and $\E$ as
the probability and expectation related to the coupled process because it is
defined as a function of $\cN$.

Let $(U_x, x\in \Z)$ be a sequence of iid random variables uniformly
distributed in $[0,1]$ and define $\sigma(x) = \one \{U_x\le \rho\}$ and
$\theta(x) = \one \{U_x\le \lambda\}$.  Call $\onu^{\lambda,\rho}$ the
resulting distribution of $(\sigma,\theta)$.  Then, $\onu^{\lambda,\rho}$ is a
product measure on $\cX^2$ with marginals
\begin{eqnarray}
  \label{p18}
  \onu^{\lambda,\rho}(\{(\sigma,\theta) : \sigma(x)=1, x\in
  A\})&=&\rho^{|A|}\nonumber\\  
  \onu^{\lambda,\rho}(\{(\sigma,\theta) : \theta(x)=1, x\in
  A\})&=&\lambda^{|A|}\nonumber 
\end{eqnarray}
for any finite $A\subset\Z$.  Furthermore, if $0\le\rho<\lambda\le 1$, then
\begin{equation}
  \label{p19}
  \onu^{\lambda,\rho}(\{(\sigma,\theta) : \sigma\le \theta \hbox{
  coordinatewise}\}) = 1 
\end{equation}

\begin{lemma} 
\label{2.1}
  The process $(\sigma_t,\theta_t)$ defined in \reff{p38} has an invariant
measure $\omu ^{\lambda,\rho}$ with marginals
\begin{eqnarray}
  \label{p22}
  \omu^{\lambda,\rho}(\{(\sigma,\theta) : \sigma(x)=1, x\in
  A\})&=&\rho^{|A|}\nonumber\\  
  \omu^{\lambda,\rho}(\{(\sigma,\theta) : \theta(x)=1, x\in
  A\})&=&\lambda^{|A|}\nonumber 
\end{eqnarray}
for any finite $A\subset\Z$.  For this measure, if $0\le\rho<\lambda\le 1$, then
\begin{equation}
  \label{pp19}
  \omu^{\lambda,\rho}(\{(\sigma,\theta) : \sigma\le \theta \hbox{
  coordinatewise}\}) = 1 
\end{equation}
Furthermore, it is possible to construct a measure $\overline\onu$ on
$\cX^2\times\cX^2$ with marginals $\omu^{\lambda,\rho}$ and
$\onu^{\lambda,\rho}$ such that
\begin{equation}
  \label{pp20}
 \overline\onu\{((\sigma,\theta),(\tsigma,\ttheta))\,:\, \sigma=\tsigma\}\,=\,1\,.
\end{equation}
\end{lemma}

\proof For the first part we follow Liggett (1985).  Start the process
$(\sigma_t,\theta_t)$ with the product measure $\onu^{\lambda,\rho}$. Since
the marginals are invariant measures for the marginal process, the first and
second marginal law of $(\sigma_t,\theta_t)$ are $\nu^\rho$ and $\nu^\lambda$,
respectively for all $t\ge0$.
% Hence we have
% $$
% \P_{\onu ^{\lambda,\rho}}(\sigma_t(x)=1\,:\, x\in A)=
% \rho^{|A|}
% $$
% $$
% \P_{\onu ^{\lambda,\rho}}(\theta_t(x)=1\,:\, x\in A)=
% \lambda^{|A|}
% $$
Therefore, we can obtain $ \omu ^{\lambda,\rho}$ as any weak limit as $t$
goes to infinity of
$$
\frac {1}{t}\int_0^t \onu^{\lambda,\rho}\bar S(s) ds.
$$
where $\bar S(t)$ is the semigroup describing the evolution of
$(\sigma_t,\theta_t)$ which is defined by $\bar\nu\bar S(t)f :=
\int\bar\nu(d(\sigma,\theta)) f(\sigma_t,\theta_t)$. The invariance of the
limit is proven in Liggett (1985).  The domination \reff{pp19} follows because
it is satisfied by the initial distribution and any transition keeps it (this
property is usually referred to as \emph{attractiveness}).

To construct a measure with the property \reff{pp20} choose $(\sigma,\theta)$
with $\omu^{\lambda,\rho}$ and define
\begin{eqnarray}
  \label{l1f}
\tsigma(x) &:=& \sigma(x)\nonumber\\
  \ttheta(x) &:=& \sigma(x)+\one\{\sigma(x)=0,\,U_x\in [\rho,\lambda]\},\nonumber
\end{eqnarray}
where $(U_x)$ is the sequence of iid uniform random variables in $[0,1]$. 
Since the $\sigma$ marginal of $\omu^{\lambda,\rho}$ is the product measure
with density $\rho$, the law of $((\sigma,\theta),(\tsigma,\ttheta))$
satisfies the requirements. \square

\section{The subadditive ergodic theorem}
Let $T:\cX^2\rightarrow {\mathcal X}$ be defined by
\[
T(\sigma,\theta)(x)=\left\{
\begin{array}{lllllllll}
1 & \hbox{if }& \sigma(x)=1  &\theta(x)=1&  &  &  &  &  \\ 
1 & \hbox{if } &  \sigma(x)=0  &\theta(x)=1& \hbox{and}& x\le 0  &  &  &  \\ 
0 & \hbox{otherwise}&  &  &  &  &  &  & 
\end{array}
\right.
\]
This operator erases the $\theta$ particles to the right of the origin, keeps
the $\theta$ particles to the left of the origin and all the $\sigma$
particles and identifies the labels $\sigma$ and $\theta$.  The operator $T$
induces a map from probability measures on $\cX^2$ to probability measures on
${\mathcal X}$ which will be called $T$ too.  Note that
\begin{equation}
  \label{pp14}
  \nu^{\lambda,\rho} = T \onu^{\lambda,\rho}
\end{equation}

%and let $\hat \mu ^{\lambda,\rho}$ br the image by $T$ of

% $\bar \mu ^{\lambda,\rho}$

\begin{propos}
\label{p14}
For all $u<v$, $\rho\leq\lambda$, there exists a random variable
$G(u,v,\lambda,\rho)$ such that:
\begin{equation}
  \label{y21}
  \P_{T\omu^{\lambda,\rho}} \left( \lim_{t\to\infty} \frac{1}{t}\sum_{ut\le x\le vt}
  \eta_t(x)=G(u,v,\lambda,\rho)\right) =1
\end{equation}
\end{propos}

The proof of this proposition follows Benassi and Fouque (1987). However, in
that paper the result is stated for the initial measure $T\onu^{\lambda,\rho}$
instead of $T\overline{\mu}^{\lambda,\rho}$.  Because the measure
$\onu^{\lambda,\rho}$ is not invariant for the process $(\sigma_t,\theta_t)$,
their random variables fail to satisfy condition~(b) of the Subaditive Ergodic
Theorem (see below).  We therefore apply that theorem to the initial measure
$T\omu^{\lambda,\rho}$ and then have to show that a strong law of large
numbers for this measure implies that a similar result holds for the more
natural initial product measure. 

We state now the subadditive ergodic theorem (taken from Liggett 1985),
introduce the notion of processes with different classes of particles and
remark some hole-particle symmetries inherent to the exclusion process. Then
we prove the proposition.

\paragraph{\bf The subadditive ergodic theorem}
{\it 
  Let $\{X_{m,n}\}$ a family of randon variables satisfying:

\noindent {\bf a.} $X_{0,0}=0$,~$X_{0,n}~\leq~X_{0,m}+X_{m,n}$~for $0\leq m 
\leq n$.

\noindent{\bf b.} $\{X_{(n-1)k,nk},~n\geq 1\}$ is a stationary sequence for 
each $k\geq 1$.

\noindent {\bf c.} $\{X_{m,m+k},~k\geq 0\}=\{X_{m+1,m+k+1},~k\geq 0\}$ in
distribution for each $m$.

\noindent {\bf d.} $\E(X_{0,1})<\infty$.

Then
$$
\lim_{n\to\infty}\frac{X_{0,n}}{n}=X_{\infty} \hbox{ exists~a.s.}
$$
}

\paragraph{\bf Particle-hole symmetry}
Holes behave as particles but with reflected rates. More precisely, given a
particle configuration $\eta$ define the \emph{reflected hole configuration}
$\check\eta$ by $\check\eta(x):=1-\eta(-x)$. The reflected hole process
$\check\eta_t$ is an exclusion process with (the same) rates $p(x,y)$. If the
particle configuration $\eta$ is distributed with $\nu^{\lambda,\rho}$, then
the reflected hole configuration $\check\eta$ is distributed with
$\nu^{1-\rho,1-\lambda}$.  Analogously, if the two-type particle configuration
$(\sigma,\theta)$ is distributed with $\mu^{\lambda,\rho}$, then the two-type
reflected hole configuration $(\check\theta,\check\sigma)$ is distributed with
$\omu^{1-\rho,1-\lambda}$, an invariant measure for the two-type reflected
hole process $(\check\theta_t,\check\sigma_t)$; furthermore
$\omu^{1-\rho,1-\lambda}$ has marginals $\nu^{1-\lambda}$ and $\nu^{1-\rho}$.
The coordinates are interchanged because $\check\theta\le\check\sigma$ and we
want to keep the first coordinate smaller than or equal to the second as in
\reff{pp19}. In particular, if $(\sigma,\theta)$ has law
$\omu^{\lambda,\rho}$, then $\check\eta=T(\check\theta,\check\sigma)$ has law
$T\omu^{1-\rho,1-\lambda}$.

\paragraph{\bf Processess with different classes of particles.}
Suppose $\sigma $ and $\xi$ are elements of $\cX$ such that $\sigma +\xi \in
\cX$, where the sum is taken coordinatewise. Then we can interpret the process
$(\eta^{\sigma}_t(\cN),\eta^{\sigma+\xi}_t(\cN)) $ as follows: a given site
$x\in Z$ is at time $t$ occupied by a first class particle, occupied by a
second class particle or vacant according to whether $\eta^{\sigma}_t(x)+
\eta^{\sigma+\xi}_t(x)$ equals $2$, $1$ or $0$ respectively. The reader can
easily check that first class particles ignore the presence of second class
particles and evolve as an exclusion process. In turn, a second class particle
jumps to empty sites as a first class particle, but when a first class
particle jumps over a second class one, they interchange sites. When $\sigma$
and $\xi$ are as above we define:
\begin{equation}
\label{p199}
\left\{
\begin{array}{llllllll}
\sigma_t=\eta^{\sigma}_t(\cN)\\
\xi_t=\eta^{\sigma+\xi}_t(\cN)-\eta^{\sigma}_t(\cN).
\end{array}
\right.
\end{equation}
So that $\sigma_t$ are the first class particles and $\xi_t$ the second class
ones.  

Let $T_m, V_m:\cX\to\cX$ be the truncations defined by
\begin{equation}
  \label{y20}
  T_m\xi(x) = \xi(x)\one\{x\le m\}\;;\qquad\qquad V_m\xi(x) =
  \xi(x)\one\{x>m\}
\end{equation}
Note that if $\sigma+\xi\in \cX$, then $T(\sigma,\sigma+\xi) =
\sigma+T_0\xi$.

\noindent{\bf Proof of Proposition \ref{p14}.}
Let $(\sigma,\theta)$
be distributed according to $\omu^{\lambda,\rho}$,
\begin{eqnarray}
\label{p201}
\xi \;:=\; T_0(\theta-\sigma)  
\end{eqnarray}
and let $(\sigma_t,\xi_t)$ be the process defined by (\ref{p199}) with initial
condition $(\sigma,\xi)$. Note that
\begin{equation}
\label{p300}
\eta_t \;:=\; \sigma_t+\xi_t,
\end{equation}
is an exclusion process with initial distribution $T\omu^{\lambda,\rho}$ and
that $\sigma_t$ is an exclusion process with initial distribution
$\nu^{\rho}$. Since this initial distribution   is  invariant,
 by the law of large numbers for triangular arrays we have
that for $n\in\N$,
\begin{equation}
  \label{26x}
 \lim_{n\to\infty}\frac1{n} \sum_{un\le x\le vn}\sigma_{n}(x) = \rho(v-u)
 \qquad\P_{\omu^{\lambda,\rho}}\hbox{ a.s.}
\end{equation}
Now taking $n=[t]$ as the integer part of $t$,
\begin{eqnarray}
  \label{y70}
  \sum_{u[t]\le x\le v[t]}\sigma_{[t]}(x) - \sum_{ut\le x\le vt}\sigma_{t}(x)
&=& \Bigl(\sum_{u[t]\le x\le v[t]}-\sum_{ut\le x\le vt}\Bigr) \sigma_{[t]}(x)
+ \sum_{ut\le x\le vt} (\sigma_{[t]}(x)-\sigma_{t}(x))\nonumber
\end{eqnarray}
The absolute value of first term is bounded by $|ut-[ut]|+|vt-[vt]|$. The
absolute value of the second term is bounded by the number of Poisson jumps
crossing either $ut$ or $vt$ in the interval $t-[t]$. Both terms converge to
zero almost surely when divided by $t$. This implies 
\begin{equation}
  \label{26}
 \lim_{t\to\infty}\frac1{t} \sum_{ut\le x\le vt}\sigma_{t}(x) = \rho(v-u)
 \qquad\P_{\omu^{\lambda,\rho}}\hbox{ a.s.}
\end{equation}
where the limit is now taken for   $t\in\R$.

In Proposition \ref{y30} below we show that for $u\ge 0$ there exists a
random variable $X(u,\lambda,\rho)$ such that
\begin{equation}
  \label{26a}
 \lim_{t\to\infty}\frac1t \sum_{x\ge ut}\xi_t(x)\;=\; X(u,\lambda,\rho)
 \qquad\P_{\omu^{\lambda,\rho}}\hbox{ a.s.}
\end{equation}
Then for $u\ge 0$ the proposition follows from \reff{p300}, \reff{26} and
  \reff{26a}; the random variable $G$ is given by
\begin{equation}
  \label{y22}
  G(u,v,\lambda,\rho) 
\;:=\; \rho(v-u) \,+\, X(u,\lambda,\rho)\,-\, X(v,\lambda,\rho)
\end{equation}

For $u<0$ we use the particle-hole symmetry. By additivity of the limits, we
can assume $v=0$. For convenience we consider $u>0$ and compute the limits for
$-u$. 
\begin{eqnarray}
  \label{y12}
  \sum_{-ut\le x< 0}\eta_t(x) 
%&=& 
% [ut] - \sum_{ut\le x< 0}
%   (1-\eta_t(x))\nonumber \\
&=& [ut] - \sum_{0< x \le ut} \check\eta_t(x)
\end{eqnarray}
where $\check\eta_t(x)$ is the reflected hole process. By the remarks made in
the particle-hole symmetry paragraph above, $\check\eta_0$ has law
$T\omu^{1-\rho,1-\lambda}$, so that we can apply the result for positive $u$
to get that \reff{y12} divided by $t$ converges $\P_{T\omu^{\lambda,\rho}}$
a.s.\/ to
\begin{equation}
  \label{y25}
 G(-u,0,\lambda,\rho)\;:=\;  u\,-\,(1-\lambda)u\,-\,
 X(0,1-\rho,1-\lambda)\,+\, X(u,1-\rho,1-\lambda)  
\end{equation}
This shows the proposition for all $u<v\in\R$. \square

\begin{propos}
  \label{y30}
  Let $u\ge 0$. Let $(\sigma,\theta)$ be  $\omu^{\lambda,\rho}$-distributed
  and let $(\sigma_t,\xi_t)$ be the two-classes particle process \reff{p199} 
with initial configuration $(\sigma,T_0(\theta-\sigma))$. 
Then there exists a random variable
  $X(u,\lambda,\rho)$ such that \reff{26a} holds.
\end{propos}

\paragraph{\bf Proof.} We prove the proposition for $u>0$. Then, the case $u=0$
follows by letting $u$ decrease to $0$ in the following inequalities:
$$ \lim_{t\to\infty}\frac1t \sum_{x\ge ut}\xi_t(x)\leq  
\lim_{t\to\infty}\frac1t \sum_{x\ge 0}\xi_t(x) \leq 
u+\lim_{t\to\infty} \frac1t \sum_{x\ge ut}\xi_t(x).$$
 We need to introduce more classes of particles.
If $\sigma $, $\xi$, $\ze$ and $\gamma$ are such that $\sigma+
\xi+\gamma+\ze\in \cX$, then
\[
(\eta^{\sigma}_t(\cN),\eta^{\sigma+\xi}_t(\cN),
  \eta^{\sigma+\xi+\gamma}_t(\cN) ,\eta^{\sigma+\xi+\gamma+\ze}_t(\cN))
\]
induces a process with first, second, third and fourth class particles by
defining, for $t\geq 0$,
\begin{equation}
\label{p200}
\left\{
\begin{array}{llllllll}
\sigma_t\,:=\,\eta^{\sigma}_t(\cN), \\ 
\xi_t\,:=\,\eta^{\sigma+\xi}_t(\cN)-\eta^{\sigma}_t(\cN), \\
\ga_t\,:=\,\eta^{\sigma+\xi+\gamma}_t(\cN)-\eta^{\sigma+\xi}_t(\cN),\\  
 \ze_t\,:=\,\eta^{\sigma+\xi+\gamma+\ze}_t(\cN)-\eta^{\sigma+\xi+\gamma}_t(\cN). 
\end{array}
\right.
\end{equation}
Let $(\sigma,\theta)$ be distributed according to $\omu^{\lambda,\rho}$ and
\begin{eqnarray}
\label{p201a}
\xi \;:=\; T_0(\theta-\sigma)   \,,\quad
\ga \;:\equiv\;  0\,,\quad
\ze \;:=\; \theta-\sigma-\xi ;
\end{eqnarray}
and let $(\sigma_t,\xi_t,\ga_t,\ze_t)$ be the four-classes process defined by
(\ref{p200}) with this initial condition. For each integer $m\geq 0$, we
define a process on $\mathcal X^4$ in the time interval $[m/u,\infty)$ as
follows: the initial (random) configurations at time $m/u$ are
 \begin{eqnarray}
\sigma^m \;:=\; \sigma_{m/u}\,,\quad
 \xi^m \;:=\; T_m(\xi_{m/u}+ \ze_{m/u} )\,,\quad
% \nonumber
% \end{eqnarray}
% \begin{eqnarray}
\ga^m \;:=\; V_m(\xi_{m/u} )\,,\quad  \ze^m \;:=\; V_m(\ze_{m/u}) \label{203}
\end{eqnarray}
and for $t\geq m/u$,
\begin{equation}
\label{p202}
\left\{
\begin{array}{llllllll}
\sigma_t^m\;:=\;\eta^{\sigma^m,m/u}_t(\cN), \\ 
\xi_t^m\;:=\;\eta^{\sigma^m+\xi^m,m/u}_t(\cN)-\eta^{\sigma^m,m/u}_t(\cN), \\
\ga_t^m\;:=\;\eta^{\sigma^m+\xi^m+\gamma^m,m/u}_t(\cN)-
\eta^{\sigma^m+\xi^m,m/u}_t(\cN), \\ 
\ze_t^m\;:=\;\eta^{\sigma^m+\xi^m+\gamma^m+\ze^m,m/u}_t(\cN)-
\eta^{\sigma^m+\xi^m+\gamma^m,m/u}_t(\cN).
 \end{array}
\right.
\end{equation}
Note that for $m=0$ this is the same process as
$(\sigma_t,\xi_t,\ga_t,\ze_t)$, and that for any $m$ the process
$(\sigma_t^m,\sigma_t^m+\xi_t^m+\ga_t^m+\ze_t^m)$ defined for $t\geq m/u$ is a
version of the coupled process \reff{p38} under the invariant measure
$\omu^{\lambda,\rho}$.

Let $u>0$ and for $0\le m\leq n$ define an  array $X_{m,n}$ of
random variables: 
\begin{equation}
X_{m,n}\;:=\;
\left\{ 
\begin{array}{llllllll}
0 & \mbox{if} & m=n \\
\sum_{y>n}\xi^m_{n/u}(y) & \mbox{if} & 0\leq m<n 
\end{array}
\right.
\end{equation} 
This array   satisfies the hypothesis of the
subadditive ergodic theorem:

\paragraph{\bf a} Assume $0<m<n$ since the other cases are trivial.
Define
\begin{equation}
  \label{y3}
\bar\xi^m \;:=\; T_m\xi_{m/u}\,.  
\end{equation}
It follows from definitions \reff{203} and \reff{y3} that for $t\ge m/u$,
\begin{eqnarray}
  \label{y1}
  \xi_t&=&\eta^{\sigma^m+\bar\xi^m+\ga^m,m/u}_t(\cN)-\eta^{\sigma^m,m/u}_t(\cN)
\nonumber\\
&\le& \eta^{\sigma^m+\xi^m+\ga^m,m/u}_t(\cN)-\eta^{\sigma^m,m/u}_t(\cN).
\nonumber
\end{eqnarray}
Hence:
\begin{eqnarray}
  \label{y2}
\xi_{n/u}&\leq&\eta^{\sigma^m+\xi^m+\ga^m,m/u}_{n/u}(\cN)-
\eta^{\sigma^m+\xi^m,m/u}_{n/u}(\cN)\;
+\;\eta^{\sigma^m+\xi^m,m/u}_{n/u}(\cN)-\eta^{\sigma^m,m/u}_{n/u}(\cN)
\nonumber\\
&=&\ga^{m}_{n/u}+\xi^{m}_{n/u}.
\nonumber
\end{eqnarray}
Therefore:
\[
X_{0,n}\;=\;\sum_{y>n}\xi_{n/u}(y)
\;\leq \;\sum_{y\in \Z}\ga^{m}_{n/u}(y)\,+\,\sum_{y>n}\xi^{m}_{n/u}(y)
\]
and, since all the $\ga$ particles are created at time $m/u$, the first term
of the right hand side  is equal to $\sum_{y\in \Z}\ga^m(y)$, which by 
(\ref{203}) is the same as $\sum_{y>m }\xi_{m/u}(y)$. Hence,
\[
X_{0,n}\leq \sum_{y>m }\xi_{m/u}(y)+\sum_{y>n}\xi^{m}_{n/u}(y)=
X_{0,m}+X_{m,n}.
\]

\paragraph{{\bf b} and {\bf c}} 
These stationary conditions follow from the space and time translation
invariance of the Poisson processes and the stationarity of the initial
measure $\omu^{\lambda,\rho}$. See the comment after \reff{p202}.

\paragraph{\bf d} 
Since $\xi$ particles can only jump when a Poisson jump is present, $X_{0,1}$,
the number of $\xi$ particles to the right of the origin at time 1, is bounded
by $\sum_{x\le0}\sum_{y> 0}(N_1(x,y)+N_1(y,x))$, the number of Poisson jumps
across the origin in the interval $[0,1]$. As in \reff{y7} and \reff{y8} (with
$k=1$) this is a Poisson variable with mean $M$ 
%$\sum_x |x| p(0,x)$ 
which is finite by hypothesis.

Then, by the subadditive ergodic theorem,
\[
\lim_{n\to\infty}\frac{X_{0,n}}{n}=X_{\infty}\quad \hbox{
  exists }\P_{\omu^{\lambda,\rho}} \hbox{ a.s.}
\]
By the definition of $X_{0,n}$, calling $n=[ut]$ (integer part)
\begin{equation}
  \label{27}
  \sum_{x\ge ut} \xi_t(x)  = X_{0,n} + \hbox{rest}
\end{equation}
where the absolute value of the rest is
\begin{eqnarray}
  \label{y4}
  \Bigl|\sum_{x\ge ut} \xi_t(x) - \sum_{x> n} \xi_{n/u}(x)\Bigr| &\le& 
\sum_{x> n} \big|\xi_t(x) -\xi_{n/u}(x)\big| \,+\,1\\
&\le& \Bigl(\sum_{x\le n}\sum_{y>n}\,+\,\sum_{y\le n}\sum_{x>n}\Bigr) (N_t(x,y)
-N_{[ut]/u}(x,y))   \,+\,1\nonumber
\end{eqnarray}
(The difference of the number of $\xi$ particles to the right of $n$ at two
different times is dominated by the number of Poisson crossings of $n$ between
those times; the addition of 1 is to cover the case $ut=n$.)  Reasoning as in
\reff{y8} (with $k=ut- [ut]$), the sum in the second line of \reff{y4} is a
Poisson random variable with mean
\[
\sum_x |x| p(0,x) (ut-[ut])/u\;\le\; \frac{M}{u}\,.
\] 
Hence, the rest divided by $t$ goes to zero
almost surely as $t\to\infty$ and we have proved that for $u>0$,
\begin{equation}
  \label{y9}
  \lim_{t\to\infty} \frac1t\,\sum_{x\ge ut}
  \xi_t(x)\,=\,X_\infty\,:=\,X(u,\lambda,\rho). \qquad\square
\end{equation}

%When $u=0$ it is possible to use the same proof defining
% \begin{eqnarray}
%\sigma^m \;:=\; \sigma_{m}\,,\quad
% \xi^m \;:=\; T_0(\xi_{m}+ \ze_{m} )\,,\quad
%\ga^m \;:=\; V_0(\xi_{m} )\,,\quad  
%\ze^m \;:=\; V_0(\ze_{m}).
%\end{eqnarray}
%and $(\sigma_t^m,\xi_t^m,\ga_t^m,\ze_t^m)$ for $t\ge m$ as
%\reff{p202} substituting $m/u$ with $m$. In this case the arrays are defined
%by 
%\begin{equation}
%X_{m,n}\;:=\;
%\left\{ 
%\begin{array}{llllllll}
%0 & \mbox{if} & m=n \\
%\sum_{y>0}\xi^m_{n}(y) & \mbox{if} & 0\leq m<n 
%\end{array}
%\right.
%\end{equation} 
%This finishes the proof of the Proposition. \square

\section{Proof of Theorem \ref{1.1}}

In the previous section we proved a law of large numbers when the initial
measure is $T\omu^{\lambda,\rho}$. In this section we show that the same is
true for $\nu^{\lambda,\rho}$ and identify the limit $G$ as being the function
$f$.

\begin{propos}
  \label{2} There exists a measure $\opi^{\lambda,\rho}$ on $\cX^2$ with marginals
  $T\omu^{\lambda,\rho}$ and $T\onu^{\lambda,\rho}$ such that the coupled
  process $(\eta_t,\teta_t)$ defined as in \reff{p38} satisfies that for all
  $K> M$,
  \begin{equation}
    \label{y40}
    \P_{\opi^{\lambda,\rho}} 
    \Bigl(\lim_{t\to\infty} \frac {1}{t} \sum_{-Kt\le x\le Kt}(\eta_t (x)-\teta_t(x))=0
  \Bigr) \; =\;1\,.
  \end{equation}
\end{propos}

\proof Choose 
$((\sigma,\theta),(\sigma,\ttheta))$ distributed according to $\overline\onu$ of
Lemma \ref{2.1} and define $\xi=\theta-\sigma$ and $\txi=\ttheta-\sigma$. Let
$\opi^{\lambda,\rho}$ be the law of 
\begin{equation}
  \label{x14}
 (\eta,\teta):=(\sigma+T_0\xi, \sigma+T_0\txi)\,.
\end{equation}
It is clear that $\opi^{\lambda,\rho}$ has marginals $T\omu^{\lambda,\rho}$ and
$T\onu^{\lambda,\rho}$. The coupling with initial distribution
$\opi^{\lambda,\rho}$  is
defined by
\begin{equation}
  \label{x5}
  (\eta_t,\teta_t)\;:=\;(\eta^\eta_t(\cN),\eta^\teta_t(\cN))
\end{equation}

Define the \emph{flux} $J^{\eta,r}_t$ as the number of $\eta$ particles that
at time zero were to the left of $r$ and at time $t$ are strictly to the right
of $r$ minus the number of $\eta$ particles that at time zero were strictly to
the right of $r$ and at time $t$ are to the left of $r$. Then for arbitrary
positive $K$ we can write
\begin{eqnarray}
  \label{y60}
\sum_{-Kt\le x\le Kt}(\eta_t (x)-\teta_t(x)) &=& \sum_{-Kt\le x\le Kt}(\eta_0
(x)-\teta_0(x)) + 
J^{\eta,-Kt}_t - J^{\eta,Kt}_t- J^{\teta,-Kt}_t +J^{\teta,Kt}_t\nonumber
\end{eqnarray}
But 
\[
\sum_{-Kt\le x\le Kt}(\eta_0
(x)-\teta_0(x))\; =\; \sum_{-Kt\le x \le0}(\theta_0
(x)-\ttheta_0(x))\,+\, \sum_{0<x\le Kt}(\sigma_0
(x)-\sigma_0(x))  
\]
By the law of large numbers for the marginals of both $T\omu^{\lambda,\rho}$
and $T\onu^{\lambda,\rho}$ the first term divided by $t$ goes to zero.  Indeed
both $\theta_0$ and $\ttheta_0$ have law $\nu^\rho$. The second term is zero.

We prove now that $\lim_{t\to\infty}(1/t)(J^{\eta,Kt}_t- J^{\teta,Kt}_t)=0$ almost
surely. 
We couple three exclusion processes with initial configurations $\sigma$, $\sigma
+ T_0\xi$ and $\sigma + T_0\txi$ and define:
\begin{equation}
  \label{y2a}
\left\{
\begin{array}{llllllll}
\sigma_t\;:=\;\eta^{\sigma}_t(\cN)\\
\xi_t\;:=\;\eta^{\sigma+T_0\xi}_t(\cN)- \eta^{\sigma}_t(\cN)\\
\txi_t\;:=\;\eta^{\sigma+T_0\txi}_t(\cN)- \eta^{\sigma}_t(\cN)\nonumber
 \end{array}
\right.
\end{equation}
In this way, $(\eta_t, \teta_t) = (\sigma_t + \xi_t, \sigma_t + \txi_t)$ and
\begin{equation}
  \label{p41}
  J^{\eta,Kt}_t- J^{\teta,Kt}_t\;=\;J^{\xi,Kt}_t- J^{\txi,Kt}_t\,.
\end{equation}
The second-class particle configurations $\xi$ and $\txi$ are dominated by
$T_0(1-\sigma)$ and dominate the null configuration. Since the dynamics is
attractive, this domination is valid at all positive times as well. This
implies that the absolute value of \reff{p41} assumes its maximal value when
$\xi= T_0(1-\sigma)$ and $\txi\equiv 0$, which we assume for the sequel. We
label the $\xi$ particles at time zero and follow the positions of the labeled
particles. Call $R^x_t$ the position of the $\xi$ particle starting at $x\le
0$. Then
\begin{equation}
  \label{a9}
  J^{\xi,Kt}_t- J^{\txi,Kt}_t 
  \;\le\;  \sum_{x\le 0} (1-\sigma(x)) \one\{R^x_t> Kt\}
\end{equation}
To dominate this consider independent Poisson random variables $N_x$, $x<0$,
with mean $1+\epsilon$. Let $Y^{x,\ell}_t$, $\ell=1,\dots,N_x$ be independent
random walks that jump from $y$ to $y+n$ at rate $p(0,n)+p(0,-n)$ for all
$y\in\Z$.  Order the $Y$ particles at time zero and call the ordered particles
$Z^i_0$ so that $Z^i_0\ge Z^{i+1}_0$ for all $i$ (here the superlabel does not
coincide necessarily with the initial position as for the $R$ and $Y$
particles). Since the mean number of $Y$ particles in each $x$ is bigger than
one, using independence and the Poisson law of $N_x$, we have $Z^i_0\ge R^i_0$
for all but a finite (random) number $W$ of $i$'s. The law of $W$ decays
exponentially.

Since the $R$ particles jump from site $x$ to site $x+n$ at most at rate 
$p(0,n)+p(0,-n)$, obvious coupling shows
\begin{equation}
  \label{a7}
  Z^i_0\ge R^i_0 \hbox{ implies } Z^i_t\ge R^i_t  \hbox{ for all } t\ge 0
\end{equation}
Then the number of $R$ particles to the right of $Kt$ is dominated by the
number of $Z$ particles to the right of $Kt$ plus $W$:
\begin{equation}
  \label{a10}
  \sum_{x\le 0}
  (1-\sigma(x))\one\{R^x_t>Kt\} \,\le\, Z_t + W
\end{equation}
where
\begin{equation}
  \label{a4}
 Z_t:= \sum_{i\le 0} \one\{Z^i_t > Kt\}\,.
\end{equation}
The variable $Z_t$ is Poisson with mean
\begin{equation}
  \label{a6}
 \E Z_t\,=\, (1+\epsilon)\sum_{x\le0} \P(Y^x_t>Kt)
\,=\, (1+\epsilon) \sum_{j\ge Kt}\P(Y^0_t>j)
%= (1+\epsilon) \sum_{j\ge 0}\P(Y^0_t>j)\one\{j\ge Kt\}
%= \E(Y^0_t\one\{Y^0_t>Kt\})
% \le (1+\epsilon) \E (Y^0_t)^\alpha
% \sum_{j\ge Kt} {1\over j^\alpha}
\end{equation}
where $Y^x_t$ is a random walk starting at $x$ that jumps from $x$ to $x+n$
with rate $p(0,n)+p(0,-n)$. Hence
% \begin{equation}
%   \label{a7}
%  \frac1t\E Z_t = (1+\epsilon) \frac1t\sum_{j\ge 0}\P(Y^0_t>j)\one\{j\ge Kt\}
% \end{equation}
\begin{eqnarray}
  \label{a13}
  \frac1t\E Z_t &=& (1+\epsilon) \frac1t\sum_{j\ge
  0}\P\Bigl(\frac{Y^0_t}t>K+\frac{j}{t}\Bigr)\\
&\le&(1+\epsilon)  \frac{[t]+1}t\sum_{j\ge
  0}\P\Bigl(\frac{Y^0_t}t>K+j\Bigr)\\
&\le&(1+\epsilon)  \frac{[t]+1}t \E\Bigl|\frac{Y^0_t}t-K\Bigr|^+
\end{eqnarray}
which goes to zero as $t\to\infty$ by the law of large numbers for the random
walk $Y^0_t$ if $K\,>\,M$, because
\begin{equation}
  \label{x15}
  M\,=\,\sum_n n (p(0,n)+p(0,-n))\,=\,\E Y^0_t/t\,.
\end{equation}

On the other hand, using the subadditive ergodic theorem we can show that
$\frac1t Z_t$ converges almost surely. The Poisson distribution of the $Y$
particles is used here to show the stationary conditions (b) and (c) of the
subadditive ergodic theorem. Indeed, the product measures with Poisson
marginals with constant mean are invariant for the process of independent
particles. Since $\frac1t Z_t$ is positive and its expectation converges to
zero, $\frac1t Z_t$ converges almost surely to zero.

This, \reff{p41}, \reff{a9} and \reff{a10}, imply
\begin{equation}
  \label{a8}
\forall K>M\qquad  \lim_{t\to\infty}
\frac1t(J^{\eta,Kt}_t-J^{\teta,Kt}_t)=0
\qquad\P_{\opi^{\lambda,\rho}}\hbox{a.s.} 
\end{equation}

The flux of particles is equal to minus the flux of holes. Indeed, each time a
particle jumps from $x$ to $y$ there is a hole jumping from $y$ to $x$. Hence,
by the particle-hole symmetry, the flux $J^{\eta,-Kt}_t-J^{\teta,-Kt}_t$ for
$(\eta,\teta)$ chosen with $\opi^{\lambda,\rho}$ has the same law as
$J^{\eta,Kt}_t-J^{\teta,Kt}_t$ for $(\eta,\teta)$ chosen with
$\opi^{1-\rho,1-\lambda}$.  This implies that
\begin{equation}
  \label{a8x}
 \forall K>M\qquad  \lim_{t\to\infty}
 \frac1t(J^{\eta,-Kt}_t-J^{\teta,-Kt}_t)=0
\qquad\P_{\opi^{\lambda,\rho}}\hbox{a.s.} \qquad \square
\end{equation}

The following is a corollary to Propositions \ref{p14} and \ref{2}.
\begin{corol}
  \label{p30}
For all $ K> M$
\[
\P_{ T \onu ^{\lambda,\rho}} {\Big (} \lim_{t\to\infty} \frac {1}{t} \sum
_{-Kt\le x\le Kt}\eta_t (x)=G(-K,K, \rho, \lambda){\Big )}=1.
\]
\end{corol}

\begin{propos}
  \label{p31}
For all $ u<v$,
% there exists a constant $H(u,v,\lambda,\rho)$ such that
%$$\lim_t \P_{  T \mu ^{\lambda,\rho}   }      {\Big (}    \vert
%\frac {1}{t} \sum _{ut}^{vt}\eta_t (x)-H(u,v, \rho, \lambda)\vert
%>\varepsilon {\Big )} 
%\rightarrow 0..$$
$$
\frac {1}{t}\sum_{ut\le x\le vt}\eta_t(x) \;\underset {t\rightarrow
  \infty}{\longrightarrow}\; \int_u^v f(s)ds\quad \hbox{ in }\ 
  \P_{T\onu^{\lambda,\rho}}\hbox{ probability }
$$
where $f$ is defined in \reff{x10}.
% $$f(u)=\left\{
% \begin{array}{lllllllll}
% \lambda & \hbox{if }& u<\alpha(1-2\lambda)  &  &  &  &  &  &  \\ 
% \frac{1}{2}(1-\frac{u}{\alpha} ) & \hbox{if } & \alpha(1-2\lambda)\leq
% u\leq   \alpha(1-2\rho) & &  &  &  &  &  \\ 
% \rho & \hbox{if}&  \alpha(1-2\rho)<u &  &  &  &  &  & 
% \end{array}
% \right. . $$
% with $$\alpha =\sum_{x\in \Z}xp(0,x).$$
\end{propos}

\paragraph{\bf Remark} 
If $p(x,y)$ vanishes when $\vert x-y\vert$ is bigger than a
constant, then the above proposition is contained in Rezakhanlou (1991).  To
include random walks with infinite range, we will derive it from Theorems
 (2.4) and (2.10)
in Andjel and Vares (1987 and 2003). 
Although their proofs are  written for the zero range
process , they  also apply to the exclusion process (see remark (5.3) in that
reference). For this process their results can be stated  as follows: 
\begin{equation}
  \label{C}
  \lim_{t\rightarrow \infty} (T \onu ^{\lambda,\rho}) S(t)\tau_{[ut]}
=   \nu^{f(u)}\ \forall u  \ \hbox{if} \  \alpha>0 \ \hbox{and }\forall u\neq
1-\lambda - \rho  \ \hbox{if}\  \alpha\leq 0, 
\end{equation}
where $f$ is given in \reff{x10} (if $\alpha >0$) and in \reff{x11} ( if $\alpha \leq 0$),
  and $\tau_x$ is the translation operator
$\tau_x\eta(y) = \eta(y-x)$. Note that 
Theorem (2.4) of Andjel and Vares (1987)
is stated for random walks with nonzero drift  but
their  proof  applies to random walks with no drift too.  The limit \reff{C} is usually referred to as
\emph{local equilibrium}. It says that at the \emph{macroscopic} point $u$ at
macroscopic time $1$ the law of the process is the invariant distribution
(equilibrium) with parameter determined by the value of the function $f$ in
site $u$.

\proof In the sequel we write $\P$ instead of $\P_{T\onu ^{\lambda,\rho}}$ to
simplify notation and $\sum_{x=a}^b$ instead of $\sum_{a\le x\le b}$ (for
$a,b\in\R$). To prove the proposition we start showing that for all $u<v$ and
$\varepsilon>0$ we have:
\begin{equation}
  \label{A}
  \lim_{t\to\infty} \P
  \Bigl[\frac {1}{t}\sum_{x=ut}^{vt}\eta_t(x)\geq
(v-u)(f(u)+\varepsilon)\Bigr]=0.
\end{equation}
Since $\frac {1}{t}\sum_{x=ut}^{(u+r)t-1}\eta_t(x)\leq r$, we may
assume that $u\neq 1-\lambda -\rho$ if $\alpha<0$ and apply \reff{C}.
Let $\delta>0$. Then, it follows from \reff{C} that there exist $n$ and $t_0$
such that:
$$
\P\Bigl[\frac
{1}{n}\sum_{x=ut}^{ut+n-1}\eta_t(x)\geq f(u)+\delta\Bigr] \leq \delta^2 \ \ \ 
\forall \ t\geq t_0.
$$
Since $(\eta_t(x), x\in \Z)$ is stochasticlly larger than $(\eta_t(x+1), x\in
\Z)$, we have 
\begin{equation}
 \label{p105}
\P\Bigl[\frac {1}{n}\sum_{x=ut+kn}^{ut+(k+1)n-1}\eta_t(x)\geq
f(u)+\delta\Bigr] \leq \delta^2 \ \ \ \forall \ t\geq t_0,k\in \N.
\end{equation}
Let
$k=k(t)=:\max\{\ell:[ut+\ell n-1\leq[vt]]\}$ then:
\begin{equation}
  \label{p33}
  \frac {1}{t}\sum_{x=ut}^{vt}\eta_t(x)=\Bigl(\frac
  {1}{t}\sum_{i=0}^{k-1}\sum_{x=ut+in}^{ut+(i+1)n-1}\eta_t(x)\Bigr)
  +\frac{g(t,\eta)}{t}\end{equation}
where $g(t,\eta)\leq n$ for all $t$.
The first term of the right hand side of \reff{p33} can be written as 
$$
\frac {nk}{t}\frac {1}{k}\sum_{i=0}^{k-1}A_i(t,\eta)
$$
where $\lim_{t\to\infty} \frac {nk}{t}=v-u$ and $A_i(t,\eta)= \frac {1}{n}
\sum_{x=ut+in}^{ut+(i+1)n-1}\eta_t(x)$.  Since
$A_i(t,\eta)\leq 1$, we have:
$$
A_i(t,\eta)\leq A_i(t,\eta)\one\{A_i(t,\eta)\leq f(u)+\delta\}+B_i(t,\eta)
$$
where $0\leq B_i(t,\eta)= \one\{A_i(t,\eta)> f(u)+\delta\}$.
  Therefore,
\begin{eqnarray}
  \label{p37}
\P\Bigl(\frac {1}{k}\sum_{i=0}^{k-1}A_i(t,\eta)\geq f(u)+2\delta\Bigr)
&\leq&
  \P\Bigl(\frac {1}{k}\sum_{i=0}^{k-1}B_i(t,\eta)\geq \delta\Bigr)\nonumber\\
&\leq& \frac
{1}{k\delta}\E\Bigl(  \sum_{i=0}^{k-1}B_i(t,\eta)  \Bigr)\;\leq \;\delta, \ \
\forall \ t\geq  t_0,
\end{eqnarray}
where the last inequality follows from \reff{p105}
Using \reff{p33} we get:
$$
\P\Bigl[\frac {1}{t}\sum_{x=ut}^{vt}\eta_t(x) \geq \frac {nk}{t} (
f(u)+2\delta)+ \frac {g(t,\eta)}{t} \Bigr]\leq \delta\ \ \forall \ t\geq t_0.
$$
Since $\lim_{t\to\infty}  \frac {nk}{t}=v-u$ and $\lim_{t\to\infty}  \frac {g(t,\eta)}{t}   =0$,
for all $t$ large enough we have:

$$
\P\Bigl[\frac {1}{t}\sum_{x=ut}^{vt}\eta_t(x) \geq (v-u) ( f(u)+3\delta)
\Bigr]\leq \delta ,
$$
 which implies \reff{A}.

Similarly one shows that 
\begin{equation}
  \label{B}
  \lim_{t\to\infty} \P\Bigl[\frac {1}{t}\sum_{x=ut}^{vt}\eta_t(x)\leq
  (v-u)(f(v)-\varepsilon)\Bigr]=0. 
\end{equation}
We now derive the proposition from \reff{A} and \reff{B}.  Let $k$ be a
positive integer and write:
$$
\frac {1}{t}\sum_{x=ut}^{vt}\eta_t(x)=\sum_{i=0}^{k-1} \frac {1}{t}
\sum_{x=ut+(v-u)t\frac {i}{k}}^{ ut+(v-u)t\frac {i+1}{k} }\eta_t(x)\,.
$$
Then apply \reff{A} to each of the terms of the sum on $i$ to obtain:
$$
\lim_{t\to\infty} \P\Bigl[ \frac {1}{t}\sum_{x=ut}^{vt}\eta_t(x) \geq
\sum_{i=0}^{k-1}\frac{v-u}{k} ( f(u+\frac {i(v-u)}{k} ) + \varepsilon)
\Bigr]=0\,,
$$
for all $\varepsilon>0$.  Letting $k$ go to infinity and then $\varepsilon$
go to $0$ we see that :
$$
\lim_{t\to\infty} \P\Bigl[ \frac {1}{t}\sum_{x=ut}^{vt}\eta_t(x) > \int_u^v
f(s)ds\Bigr]=0\,.
$$
Similarly, using \reff{B} we get
$$
\lim_{t\to\infty} \P\Bigl[ \frac {1}{t}\sum_{x=ut}^{vt}\eta_t(x) < \int_u^v
f(s)ds\Bigr]=0\,,
$$
and the proposition follows. \square

\begin{corol}
  \label{5}
For all $ K> M$
$$
\P_{ T \onu^{\lambda,\rho} } {\Big (} \lim_{t\to\infty} \frac {1}{t} \sum
_{-Kt\le x \le Kt}\eta_t (x)=\int_{-K}^K f(s)ds{\Big )}=1\,.
$$
\end{corol}

\begin{propos}
  \label{6}
For all $ u<v$ and $\rho \leq \lambda$
$$
\P_{ T\onu^{\lambda,\rho} } {\Big (} \liminf_t \frac {1}{t} \sum
_{ut\le x \le vt}\eta_t (x)\geq \int_u^v f(s)ds {\Big )}=1\,.
$$
\end{propos}

\proof 
%We will first prove the proposition for $\alpha (1-\rho )\leq u<v:$
Let $(\sigma_0,\theta_0)$ be distributed according to
 $\onu^{\lambda,\rho}$ . Then, let
 $\xi_0 = T_0(\theta_0-\sigma_0)$ and let 
 $(\sigma_t,\xi_t)$ be the process defined by (\ref{p199}) with this initial
 condition. Then,
\begin{equation}
\label{p400}
\eta_t \;:=\; \sigma_t+\xi_t
\end{equation}
is the exclusion process with initial distribution $T\onu^{\lambda,\rho}$, and
$\sigma_t$ is the exclusion process with initial distribution
$\nu^{\rho}$. Then,  as in the proof of \reff{26},
$$
\P_{\onu^{\lambda,\rho}}\Bigl(\lim_{t\to\infty} \frac{1}{t} \sum_{ut\le
  x \le vt} \sigma_{t}(x)=\rho (v-u)\Bigr)=1. 
$$
This and \reff{p400} imply that
$$
\P_{ T \onu ^{\lambda,\rho} } {\Big (} \liminf_{t\to\infty} \frac {1}{t} \sum _{ut\le
  x \le vt}\eta_t (x)\geq \rho (v-u) {\Big )}=1\,.
$$
Since $f(s)=\rho$ if $s\geq \alpha(1-2\rho)$, this proves the proposition
for $\alpha (1-2\rho )\leq u<v$.

For the case $u<v\leq \alpha (1-2\rho )$, note that 
$T \onu ^{\lambda,\rho}\geq T\onu^{\lambda,0}=T\omu^{\lambda,0}$. Thus, the 
result  for these values of $u$ and $v$ follows
from Proposition \ref{p31} and the fact that the function  $f$ in
\ref{x10} does not change its values in the interval
 $(-\infty,\alpha (1-2\rho )]$ if we substitute $0$ for $\rho$

Finally for $u<\alpha (1-2\rho )<v$ the result follows from the inequality
$$\liminf_{t\to\infty} \frac {1}{t} \sum_{x=ut}^{vt}\eta_t (x)\geq
\liminf_{t\to\infty} \frac {1}{t} \sum_{x=ut}^{\alpha (1-2\rho )t}\eta_t (x)+
\liminf_{t\to\infty} \frac {1}{t} \sum_{x=\alpha (1-2\rho )t}^{vt}\eta_t (x)$$
 and the two previous cases. \square

\bigskip

\paragraph{\bf Proof of Theorem \ref{1.1}} 
Fix $u<v$ and let $K= 1+\max\{M,\vert u \vert , \vert v \vert \}$ Then, by
Corollary \ref{5} we have $\P_{ T \onu ^{\lambda,\rho} }$-a.s.
\begin{eqnarray}
  \label{w1}
\int_{-K}^{K} f(s)ds&=&\lim_{t\to\infty} \frac {1}{t} 
\sum_{x=-Kt}^{Kt}\eta_t (x)\nonumber\\
&\geq&
\liminf_{t\to\infty} \frac {1}{t} \sum _{x=-Kt}^{ut-1}\eta_t (x)+
\limsup_{t\to\infty} \sum _{x=ut}^{vt}\eta_t (x) + \liminf_{t\to\infty}
 \frac {1}{t} \sum
_{x=vt+1}^{Kt}\eta_t (x ) \nonumber
\end{eqnarray}
Therefore, by Proposition \ref{6}
$$
\limsup_{t\to\infty} \sum _{x=ut}^{vt}\eta_t (x)\leq \int_{-K}^K f(s)ds- \int_{-K}^u
f(s)ds- \int_{v}^K f(s)ds= \int_{u}^v f(s)ds \ \
P_{ T \onu ^{\lambda,\rho} }- a.s.
$$
The theorem follows from Proposition \ref{6}, this last inequality and 
(\ref{pp14}).
\square

\section*{Acknowledgements} 
We thank Maria Eulalia Vares for useful discussions concerning her paper with
the first author. This work is partially supported by the Brazilian agencies
FAPESP, CNPq and FINEP, the Brasil-France Mathematical Network and the
USP-COFECUB agreement.

\vskip 1cm
\parindent 0pt

%\address{$^a$
Enrique D. Andjel\newline
CMI, Universit\'e de Provence\newline 
39 rue Joliot Curie\newline 
13453 Marseille cedex 13\newline 
FRANCE\newline
\texttt{Enrique.Andjel@cmi.univ-mrs.fr}
%}
\vskip 4mm

%\address{$^b$
Pablo A. Ferrari and Adriano Siqueira \newline
IME USP,\newline 
Caixa Postal 66281, \newline 
05311-970 - S\~{a}o Paulo, BRASIL. \newline 
\rm \texttt{pablo@ime.usp.br}, \texttt{adr@ime.usp.br},\newline 
http://www.ime.usp.br/\~{}pablo
%}

 %%%%%%%%%%%%%%%%%%%%%%%%%%%%%%%%%%%%%%%%%%%%%%%%%%%%%%%%%%%%%
\end{document}